\documentclass[runningheads,a4paper]{llncs}

\usepackage{hyperref}
\usepackage{amssymb}
\usepackage{graphicx}
\usepackage{amsmath}
\usepackage{verbatim}
\usepackage{amssymb}
\usepackage[sectionbib]{natbib}
\usepackage{url}
\usepackage{xr}
\usepackage{fullpage}
\usepackage{appendix}
\usepackage{color}
\usepackage{float}

\usepackage{amsthm}
\newtheorem*{theoremNN}{Theorem}

\makeatletter
\def\@seccntformat#1{\@ifundefined{#1@cntformat}%
   {\csname the#1\endcsname\quad}  
   {\csname #1@cntformat\endcsname}
}
\let\oldappendix\appendix 
\renewcommand\appendix{%
    \oldappendix
    \newcommand{\section@cntformat}{\appendixname~\thesection\quad}
}
\makeatother

\urldef{\mailsa}\path|{hamedkarim@gmail.com},{jnutini, schmidtm}@cs.ubc.ca|

\def\norm#1{\|#1\|}
\newcommand{\R}{\mathbb{R}}

\newcommand{\half}{\frac 1 2}

\newcommand{\be}{\begin{equation}}
\newcommand{\ee}{\end{equation}}
\newcommand{\bea}{\begin{eqnarray}}
\newcommand{\eea}{\end{eqnarray}}
\newcommand{\ben}{\begin{equation*}}
\newcommand{\een}{\end{equation*}}
\newcommand{\bean}{\begin{eqnarray*}}
\newcommand{\eean}{\end{eqnarray*}}

\definecolor{red}{rgb}{1,0,0}

\newcommand{\ru}[1]{\rotatebox[origin=c]{25}{$#1$}}   
\newcommand{\rd}[1]{\rotatebox[origin=c]{-25}{$#1$}}  

\begin{document}

\mainmatter  

\title{The Polyak-\L{}ojasiewicz Paper: A 10-Year Retrospective}


%
%
\author{Hamed Karimi,  Julie Nutini, Mark Schmidt%
}
%

\institute{Department of Computer Science, University of British Columbia (as of 2016)\\
}

%
%

\maketitle

\begin{abstract}
In 2016, we wrote a paper introducing the phrase ``Polyak-\L{}ojasiewicz (PL) condition'' for an old idea, unsure of how readers would react.
The paper ultimately had a larger impact than we expected, and we are surprised and honoured to receive the 2026 PKDD ECML Test of Time Award.
With 10 years of hindsight, this note offers an informal and candid retrospective on the PL condition. 
We take a look back at the origins of the paper, give an AI-assisted refined version of the main theorem, review criticisms of the work, and discuss a variety of related ideas and extensions.
\end{abstract}

\section{Story of the Paper}\label{sect:intro}

\begin{figure}
\begin{center}
\includegraphics[width=.7\textwidth]{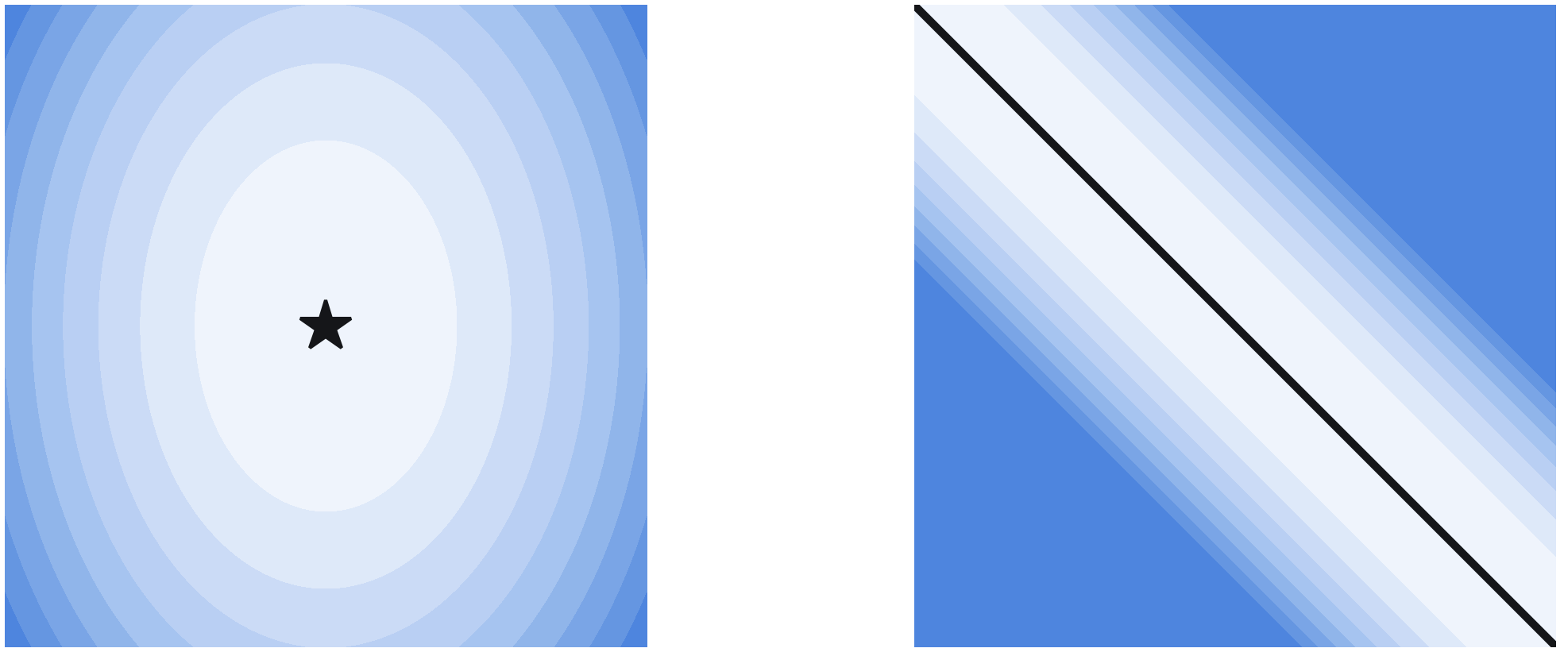}
\caption{A two-dimensional strongly-convex objective (left) and a least squares objective with 2 features (right). Even with 2 features, least squares may have an infinite number of solutions and may not satisfy strong convexity.}
\label{fig:LS}
\end{center}
\end{figure}

Our work on the PL paper started with a simple question: 
\begin{quote}``\emph{Which assumption gives the simplest proof showing that gradient descent (GD) has a linear convergence rate for fitting linear least squares models?}''\end{quote}
As context, we usually fit machine learning (ML) models with gradient descent (or one of its variations like stochastic gradient descent or Adam). The textbook result about the speed of gradient descent is that it has a linear convergence rate for functions that are strongly-convex. However, strong convexity does not hold for most ML models. For example, as illustrated in Figure~\ref{fig:LS}, strong convexity does not always hold for the ``least squares'' problem of fitting a linear regression model with the squared error. 

\vspace{-2em}

The standard alternative to strong convexity was the \emph{error bound} condition of~\citet{Luo}, but there also existed the \emph{quadratic growth} condition~\citep{bonnans1995second} and a variety of newer alternatives to show that gradient descent can achieve a linear convergence rate:
\begin{quote}
\emph{Restricted secant inequality} and \emph{restricted strong convexity}~\citep{zhang-yin}, \emph{semi-strong convexity}~\citep{gong2014linear}, \emph{essential strong convexity}~\citep{Liu_esc}, \emph{optimal strong convexity}~\citep{Liu_osc}, \emph{weak strong convexity} under two different definitions~\citep{Necoara2015, Ma2015}, and later \emph{quasar strong convexity}~\citep{hinder2020near}, \emph{strong star convexity}~\citep{zhang2023prise}, \emph{strong quasi convexity}~\citet{lara2025characterizations} which was originally studied by~\citet{polyak1966existence} in a different context, and \emph{strong star quasi convexity}~\citep{khanh2025star}.
\end{quote}
This variety of conditions led us to the question above. But the answer ended up being that there was a condition that gave a simpler proof than all of the above assumptions. And this condition was weaker than all of the above assumptions for giving linear convergence. And this condition did not even imply convexity.

\subsection{Polyak-\L{}ojasiewicz Condition}

\begin{figure}
\begin{center}
\includegraphics[width=.26\textwidth]{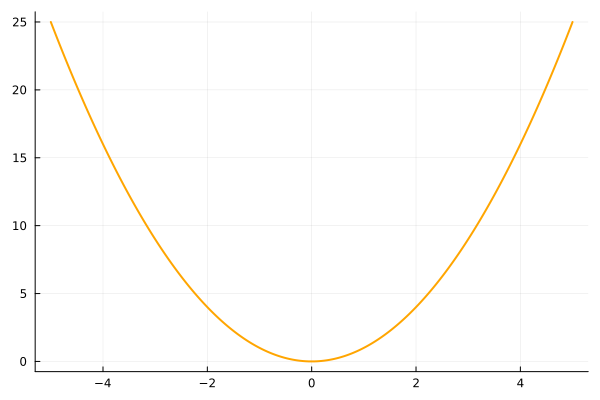}
\includegraphics[width=.26\textwidth]{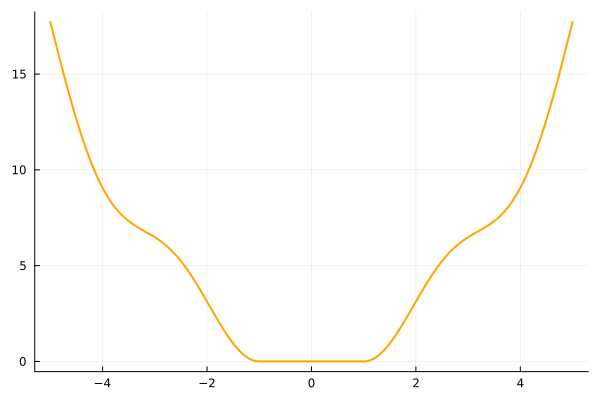}
\includegraphics[width=.26\textwidth]{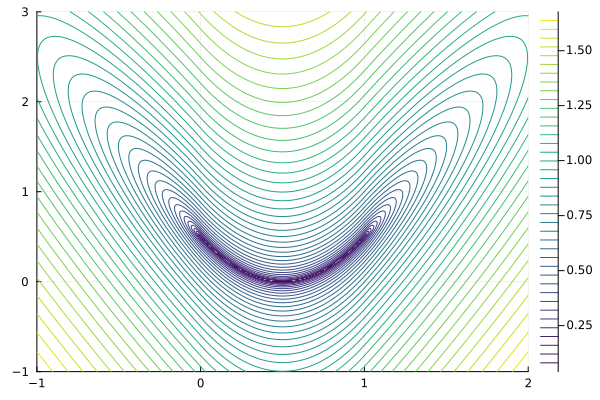}
\caption{A  strongly-convex function (left) and two non-convex PL functions (middle and right). The left function is $x \to x^2$ (PL with $\mu=2$). The middle function is $x \to d(x)^2 + 3\sin^2(d(x))$ with $d(x) = \max\{0, |x| - 1\}$ which has an interval of minimizers (PL with $\mu=1/32$). The third is $(x,y) \to d(x)^2 + (y - 2(x - 1/2)^2 + 2d(x)^2)^2$ with $d(x)=\max\{0, x-1, -x\}$ from Gabriel Peyré and has a curved valley of solutions (PL with $\mu=6-4\sqrt{2}$).  An interactive animated variant on the right plot is available here:\\ \url{https://www.gpeyre.com/mathematical-nexus/python/p-l-flow/banana_gradient_flow_interactive.html}}.
\label{fig:examples}
\end{center}
\end{figure}

The specific problem we focused on is minimizing a differentiable function $f$ that is bounded below by some value $f^*$ and that has parameters $x \in \R^n$ . The specific condition we considered was that the norm of the gradient of $f$ grows at least quadratically in the sub-optimality,
\be
\label{pl_ineq}
\frac{1}{2}|| \nabla f(x) ||^2 \geq \mu (f(x) - f^*), \quad \forall~x \in \R^n.
\ee
We require that this inequality holds for a positive $\mu$ for all $x$. This inequality can be used in many proofs at the point where you would normally use strong convexity. However, unlike strong convexity, the above inequality always holds for least squares, allows multiple solutions, and allows the function to be non-convex (see Figure~\ref{fig:examples}). We thought this was interesting enough to write up as a short paper. 

\begin{figure}
\begin{center}
\includegraphics[width=.4\textwidth]{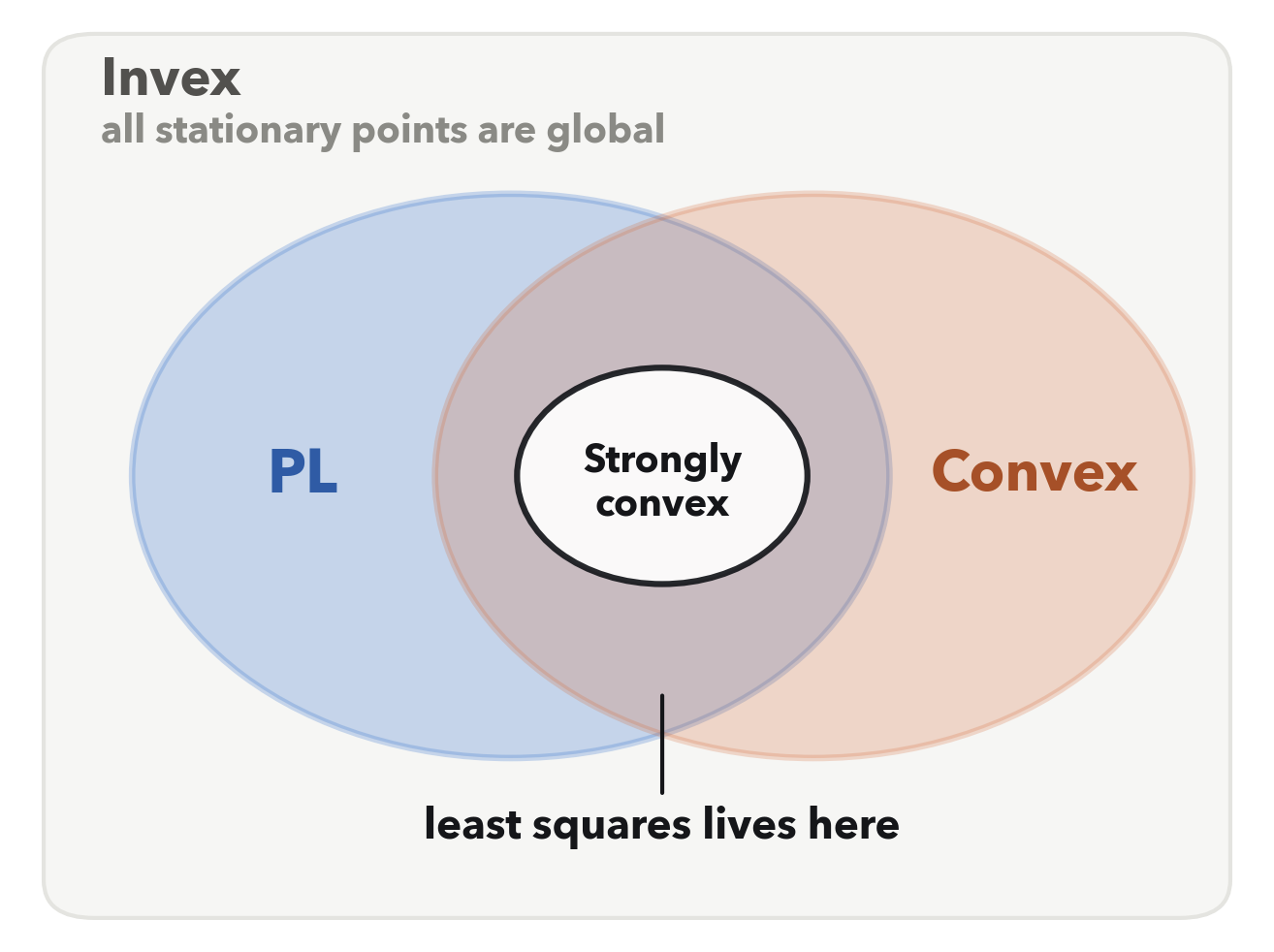}
\caption{The PL assumption implies invexity but not convexity. Convexity also implies invexity. Strong-convexity implies both PL and convexity. Linear least squares problems are at the intersection of PL and convexity. }
\label{fig:invex}
\end{center}
\end{figure}

While the condition in~\eqref{pl_ineq} does not imply convexity, it does imply invexity which for smooth functions means that all stationary points (where $\nabla f(x)=0$) are global optima (see Figure~\ref{fig:invex}). Thus, we were initially calling it the ``strong invexity'' condition. However, Stephen Wright (who referred to the property as the ``strong property'' in his ISMP 2015 lecture) pointed us to the invexity literature where we found that ``strong invexity'' was already taken. As we continued to write our paper, we soon found that the condition was already published twice in 1963, in Russian by~\citet{polyak} and in French by~\citet{loj}. Polyak gave the condition as an example of an assumption that yields linear convergence for gradient descent. 
 \L{}ojasiewicz showed that any real-analytic function satisfied such a condition locally near any critical point (possibly with a different exponent), and their work is viewed as foundational in many fields. 
As~\ref{pl_ineq} was a special case of the condition referred to as the Kurdyka-\L{}ojasiewicz (KL) condition studied by~\citet{bolte2010characterizations}, we chose to follow their naming convention and referred to it as the Polyak-\L{}ojasiewicz (PL) condition. While preparing this note, J\'{e}rome Bolte highlighted that~\citet{lezanski1963minimumproblem} also discusses (in German) such a condition to show linear convergence of gradient methods in Banach spaces.

\subsection{Writing/Submitting/Resubmitting}

The paper~\citep{karimi2016linear} was difficult to write because the contribution was unconventional. The PL condition was already known in the 1960s, and these ideas seemed reasonably well explored for convex functions and for iterates near the solution (even in the non-smooth case). In an attempt to make future reviewers happy, we chose to focus on showing that the PL condition was the weakest condition in the (smooth) non-convex setting. We also came up with the proximal-PL condition to simplify proofs in the proximal-gradient setting, and added PL-replacing-strong-convexity proofs for many standard algorithms used in ML. We first presented the work at the 2015 NeurIPS Optimization for ML workshop, followed by the CIFAR Neural Computation and Adaptive Perception program meeting. The overall feedback was that people liked knowing that this was the weakest condition and that it could easily be swapped in for strong convexity. However, the reception of the work seemed fairly lukewarm at these venues, and there were some early criticisms (expanded on in Section~\ref{sec:crit}).

Beyond the workshop paper, we submitted the work to AISTATS where it was rejected. Below are some quotes from the mixed reviews:
\begin{quote}
{\bf Reviewer 1}: “\emph{whether it is really useful in practice is unclear to me}.”
\end{quote}
\begin{quote}
{\bf Reviewer 2}: “\emph{paper should be beneficial to this field and be publicly known}.”
\end{quote}
\begin{quote}
{\bf Reviewer 3}: “\emph{from the technical view point, the contribution is limited in a minor modification}.”
\end{quote}
In light of these reviews and a technical error that came up during the discussion (our argument that a particular problem class satisfied the assumption was incorrect), the paper was rejected. After fixing our error, we submitted the paper to ECML. While the reviews were still mixed, they were shifted slightly more favourably:
\begin{quote}
{\bf Reviewer 1}: “\emph{interesting although has limited novelty as it mostly builds on previous work}.”
\end{quote}
\begin{quote}
{\bf Reviewer 2}: “\emph{realization of the different kind of learning problems to achieve better performance}.”
\end{quote}
\begin{quote}
{\bf Reviewer 3}: “\emph{this paper can bring a large impact on machine learning and optimization community}.”
\end{quote}
The paper was ultimately accepted to ECML. After we put the paper on arXiv for the conference, it seemed like we immediately started hearing the expression ``PL condition'' in a variety of contexts. We did not expect the paper to have much impact. We certainly did not expect it to accumulate several thousand citations and start appearing in textbooks. We were also surprised to receive this message from Polyak himself around two weeks after we put the first draft on arXiv:

\begin{center}
\includegraphics[width=.9\textwidth]{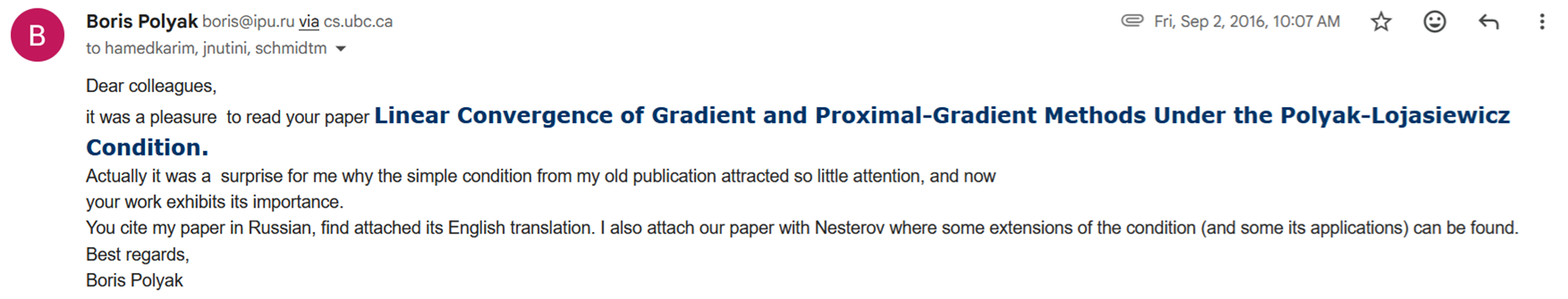}
\end{center}

\subsection{Contributions}

The final paper ultimately contained 3 main contributions.
\begin{enumerate}
\item We gave the relationship between conditions for the case of smooth and not-necessarily-convex functions, showing that the old PL condition and the equivalent-but-also-relatively-old EB condition are the most general conditions under which gradient descent has a linear convergence rate. We review these relationships in Section~\ref{sec:thm2}, and give a refined version of the earlier result. Most of these relationships were already in the literature for the special case of convex functions, and all of these relationships except one have a straightforward 2-line proof. The difficult proof (showing that either PL or EB implies quadratic growth) was shown in 2015 for convex functions in two works~\citep{Bolte2015,Zhang2015}. Showing this result in the non-convex setting was done independently of our work by~\citet{zhang2016characterization}, and the construction used is almost identical between our works (we ultimately had to check our records to verify that we had the construction before their paper appeared on arXiv). On the other hand,~\citet{rebjock2025fast} argue that PL/EB implying QG appears earlier in the literature, including as early as 2000 in works by Ioffe and by Otto \& Villani.
\item We gave a non-smooth generalization of the PL condition for analyzing proximal-gradient methods called the ``proximal-PL condition". Proximal-gradient methods apply to functions of the form $F = f + g$, where $f$ is smooth and $g$ is convex. We showed that the proximal-PL condition was equivalent to the older EB and KL conditions in the proximal-gradient setting, but the proximal-PL condition yields a simpler proof. However, showing that a function satisfies the proximal-PL condition may not be easy at all; we gained a simple proof by giving up a simple assumption. Nevertheless, we gave 5 problem classes where the proximal-PL assumption is guaranteed to be satisfied. Although these include many interesting cases such as L1-regularized least squares, they are less general than one might hope. For example, while $f$ being strongly-convex guarantees we satisfy proximal-PL, $F$ being strongly convex does not guarantee that the proximal-PL assumption is satisfied. There are even counter-examples showing that the proximal-gradient method does not obtain a linear convergence rate in general in this case.
\item We analyzed common algorithms in ML under the PL and proximal-PL condition. This included coordinate descent (both randomized and greedy versions), SignGD (a special case of the Adam optimizer), stochastic gradient descent (under decreasing step sizes, constant step sizes, and variance reduction), and proximal coordinate descent (which includes dual solvers for support vector machines). In all cases, we obtained the strongly-convex convergence rate under the weaker assumption. At the time we thought that you could replace strong convexity with PL in the proofs for every algorithm. While this ended up being true in many cases, it turned out to not be true in some important cases like accelerated gradient methods (we discuss this further in Sections~\ref{sec:crit} and~\ref{sec:ext}).
\end{enumerate}

\section{Relationships Between Conditions (Again)}
\label{sec:thm2}

Below is the result that the paper~\citep{karimi2016linear} gives showing  relationships between the conditions.
\begin{theoremNN}
For a function $f$ with a Lipschitz-continuous gradient, the following implications hold:
\[
(SC) \subseteq (ESC) \subseteq (WSC) \subseteq (RSI) \subseteq (EB) \equiv (PL) \subseteq (QG).
\]
If we further assume that $f$ is convex, then we have 
\[
(RSI) \equiv (EB) \equiv (PL) \equiv (QG).
\]
\end{theoremNN}
In the statement, (SC) refers to the class of functions with a Lipschitz-continuous gradient that also satisfy strong-convexity
and similarly for the other abbreviations. 
The paper also discusses different names for these conditions, and notes that PL implies invexity but QG does not. Note that QG is the only condition that does not yield linear convergence to a global minimizer of gradient descent.

Notice that many of the conditions are equivalent for convex functions, while EB and PL are unconditionally equivalent given the Lipschitz assumption. In the paper we use $x \to x^2 + 3\sin^2(x)$ as an example of a PL function that is non-convex, but this function also satisfies the RSI condition. This is related to a question we were curious about in 2016 but did not fully explore: are there other equivalences between these classes?

In this section we first define the various conditions and over a dozen names they have appeared under in the literature. Then we give a refined version of the above result. In particular, we give a sequence of simple counter-examples showing that there are no other equivalences in the general case, even if we assume invexity. We also add 3 conditions to the hierarchy (quasar strong convexity, strong quasi convexity, and strong star quasi convexity) and show that the hierarchy collapses further  in the convex setting (and even further for convex quadratic functions). We note that Gemini was used to help construct these arguments, an option that was not available in 2016.

\subsection{[Multiple-]Definitions of the Conditions}

Here we give the definitions we use for the various conditions, largely copied from Appendix A of the original paper. We also mention alternate names they have been used under in the literature. We focus on the case of differentiable functions, although these conditions have non-differentiable generalizations. Note that all conditions require the constant $\mu$ to be positive, and that the constant $\mu$ \emph{is not necessarily the same} across different conditions. 

Many of the conditions require a reference point $x_p$ that is in the solution set. Although other definitions exist in the literature, we define $x_p$ as the projection onto the solution set. In particular, we define $x_p$ as the value $w$ minimizing $\norm{x-w}$ among the set $\arg\min_z\{f(z)\}$. An issue with this definition of $x_p$ for non-convex functions is that this projection may not be unique. Although other ways to deal with this are possible, we will assume that the condition holds for \emph{every nearest minimizer} (all $w$ minimizing $\norm{x-w}$ in the argmin). Some papers require the inequality to hold for all $x_p$ in the solution set~\citep[see][]{khanh2025star}, but this is a stronger assumption and for many conditions would require that the solution is unique.

\begin{enumerate}
\item \textbf{Strong Convexity} (SC): For all $x$ and $y$, we have
\[
f(y) \geq f(x) + \langle \nabla f(x),y-x\rangle + \frac \mu 2 \norm{y-x}^2.
\]
For $C^2$ functions, SC is equivalent to $\nabla^2 f(x)$ having eigenvalues at least $\mu$ for all $x$. Note that SC implies convexity. \smallskip
\item \textbf{Essential Strong Convexity} (ESC): For all $x$ and $y$ with the same projection onto the solution set ($x_p = y_p$) we have~\citep{Liu_esc}
\[
f(y) \geq f(x) + \langle \nabla f(x),y-x\rangle + \frac \mu 2 \norm{y-x}^2,
\]
Note that ESC, similar to the conditions below, does not imply convexity (see Appendix~\ref{sec:convex}). \smallskip
\item \textbf{Weak Strong Convexity} (WSC): For all $x$, we have~\citep{Necoara2015}
\[
f^* \geq f(x) + \langle \nabla f(x),x_p-x\rangle + \frac \mu 2 \norm{x_p-x}^2.
\]
This has also been called \textbf{strong star convexity}~\citep{zhang2023prise} and \textbf{quasi-strong convexity}~\citep{hui}. \smallskip
\item \textbf{Strong Quasi Convexity} (SQC): For all $x$ and $y$, we have for all $\lambda\in[0,1]$ that~\citet{polyak1966existence}
\[
f(\lambda x + (1-\lambda)y) \leq \max\{f(x), f(y)\} - \frac \mu 2\lambda(1-\lambda)\norm{y-x}^2.
\]
For differentiable functions, this is equivalent to
\[
0 \geq  \langle \nabla f(x),y -x \rangle + \frac \mu 2 \norm{y-x}^2,
\]
for all $f(y) \leq f(x)$. This equivalence was first stated in one of Yuri Nesterov's earliest publications~\citep{VNCh_1978}, although the derivation there has an issue that was amended by~\citet{hadjisavvas2026characterizations}.
Note that SC and SQC are the only conditions implying that the solution is unique.
\item \textbf{Strong Star Quasi Convexity} (SSQC): For all $x$, we have for all $\lambda\in[0,1]$ that~\citet{khanh2025star}
\[
f(\lambda x + (1-\lambda)x_p) \leq f(x) - \frac \mu 2\lambda(1-\lambda)\norm{x-x_p}^2.
\]
\item \textbf{Quasar Strong Convexity} (QSC): For all $x$ and all $\lambda \in [0,1]$, we have~\citep{hinder2020near}
\[
f(\lambda x_p + (1-\lambda)x) \leq \gamma\lambda f^* + (1-\gamma\lambda)f(x) - \frac{\gamma\mu}{2}\lambda\left(1-\frac{\lambda}{2-\gamma}\right)\norm{x_p-x}^2,
\]
for some $\gamma \in (0,1]$. Or equivalently, for differentiable $f$ we have
\[
f^* \geq f(x) + \frac 1 \gamma \langle \nabla f(x),x_p-x\rangle + \frac \mu 2 \norm{x_p-x}^2. \smallskip
\]
\item \textbf{Restricted Secant Inequality} (RSI): For all $x$, we have~\citep{zhang-yin}
\[
\langle \nabla f(x),x-x_p\rangle \geq \mu \norm{x_p-x}^2.
\]
This is also called the \textbf{correlated-point error bound}~\citep{zhang2016characterization}.
If the function $f$ is also convex, then it is called \textbf{restricted strong convexity}~\citep{zhang-yin}. \smallskip
\item \textbf{Error Bound} (EB): For all $x$, we have
\[
\norm{\nabla f(x)} \geq \mu\norm{x_p -x}.
\]
This is also called the \textbf{residual-point error bound}~\citep{zhang2016characterization}. \smallskip
\item \textbf{Polyak-\L{}ojasiewicz} (PL): For all $x$, we have
\[
\half \norm{\nabla f(x)}^2 \geq \mu(f(x) - f^*).
\]
This is a special case of being \textbf{gradient dominated}~\citep{nesterov2006cubic}, and is called the \textbf{residual-objective error bound}~\citep{zhang2016characterization}. Notice that this condition does not require dealing with the issue of defining $x_p$, which is appealing for non-convex problems. \smallskip
\item \textbf{Quadratic Growth} (QG): For all $x$, we have~\citep{bonnans1995second}
\[
f(x) - f^* \geq \frac \mu 2 \norm{x_p-x}^2.
\]
This is called the \textbf{objective-point error bound}~\citep{zhang2016characterization}. If the function $f$ is also convex, then this condition has also been called \textbf{semi-strong convexity}~\citep{gong2014linear}, \textbf{optimal strong convexity}~\citep{Liu_osc}, and in one work is called weak strong convexity~\citep{Ma2015}.
\end{enumerate}
We note that up to a re-scaling of $\mu$ the conditions can summarized as:
\begin{center}
\begin{tabular}{cclcccl}
(SC) & $\langle \nabla f(x), x - y\rangle $ & $\geq$ & $ f(x)  - f(y)$  & + & $\frac \mu 2 \norm{x-y}^2$ & for all $x$ and $y$.\\
(ESC) & $\langle \nabla f(x), x - y\rangle $ & $\geq$ & $ f(x)  - f(y)$   &+ & $\frac \mu 2 \norm{x-y}^2$ & for all $x$ and $y$  with $x_p=y_p$.\\
(SQC) & $\langle \nabla f(x), x - y\rangle $ & $\geq$ & 0 & + &  $\frac \mu 2 \norm{x-y}^2$ & for all $x$ and $y$ with $f(y) \leq f(x)$.\\
(WSC) & $\langle \nabla f(x), x - x_p\rangle $ & $\geq$ & $ f(x)  - f^*$ &  + & $\frac \mu 2 \norm{x-x_p}^2$ & for all $x$.\\
(QSC) & $\langle \nabla f(x), x - x_p\rangle $ & $\geq$ & $ \gamma(f(x)  - f^*)$ &  + & $\frac \mu 2 \norm{x-x_p}^2$ & for all $x$ and some $\gamma \in [0,1]$\\
(RSI) & $\langle \nabla f(x), x - x_p\rangle $ & $\geq$ & 0  & +  & $\frac \mu 2 \norm{x-x_p}^2$ & for all $x$\\
(EB) & $\langle \nabla f(x), \nabla f(x)\rangle $  & $\geq$ & 0 & + & $\frac \mu 2 \norm{x-x_p}^2$ & for all $x$\\
(PL) & $\langle \nabla f(x), \nabla f(x)\rangle $ & $\geq$ & $  \mu(f(x) - f^*)$  & + & 0 & for all $x$\\
(QG) & $f(x)-f^*$  & $\geq$ & 0 & + & $\frac \mu 2 \norm{x-x_p}^2$ & for all $x$\\
\end{tabular}
\end{center}
We omit SSQC above which we will see is equivalent to 2 other conditions for smooth functions. From this summary we see that all conditions are equivalent to a lower bound on particular directional derivatives, except QG which only lower bounds function values.

\begin{figure}[H]
\begin{center}
\includegraphics[width=.32\textwidth]{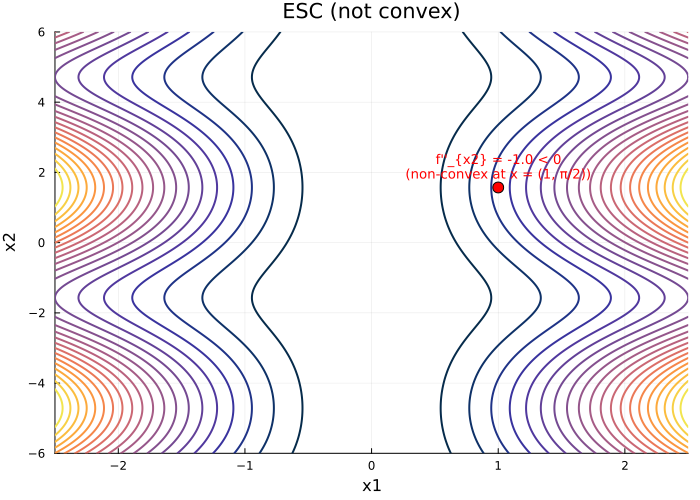}
\includegraphics[width=.32\textwidth]{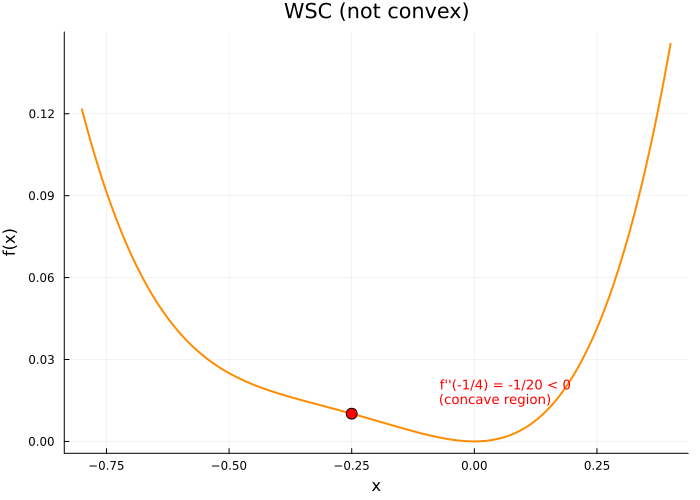}
\includegraphics[width=.32\textwidth]{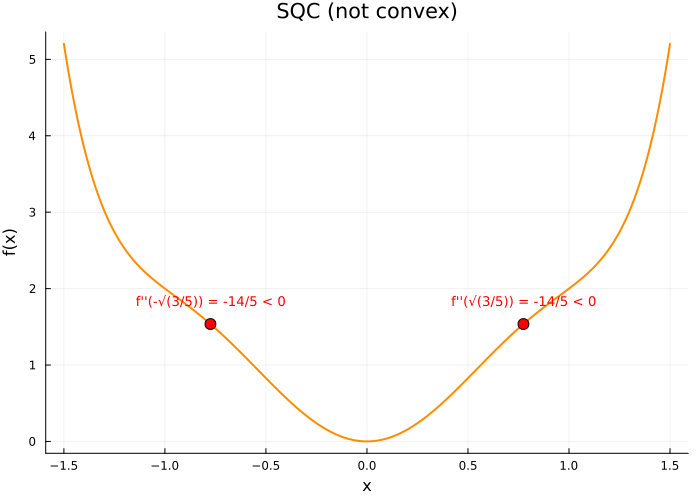}
\caption{Non-convex functions satisfying the strongest conditions besides SC. ESC functions can only be non-convex in dimension 2 or higher.}
\label{fig:convex}
\end{center}
\end{figure}

We note that among these 10 conditions, \emph{only SC implies convexity}, except in one dimension where ESC also implies convexity. We show  this (assisted by Gemini) in Appendix~\ref{sec:convex} for the strongest conditions. Figure~\ref{fig:convex} shows non-convex counter-examples for the strongest conditions.

\subsection{Updated Result}
\label{sec:updated}

\begin{figure}[H]
\begin{center}
\includegraphics[width=.4\textwidth]{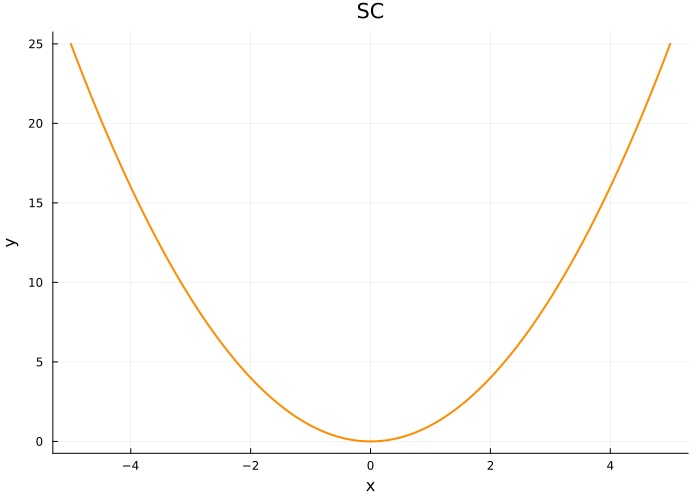}
\includegraphics[width=.4\textwidth]{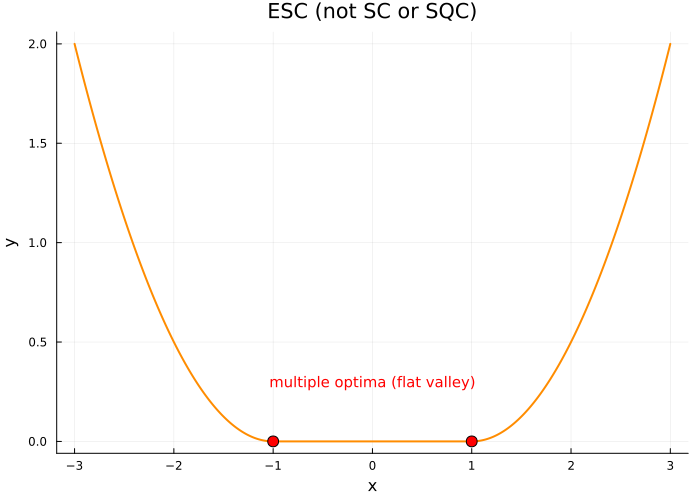}
\includegraphics[width=.4\textwidth]{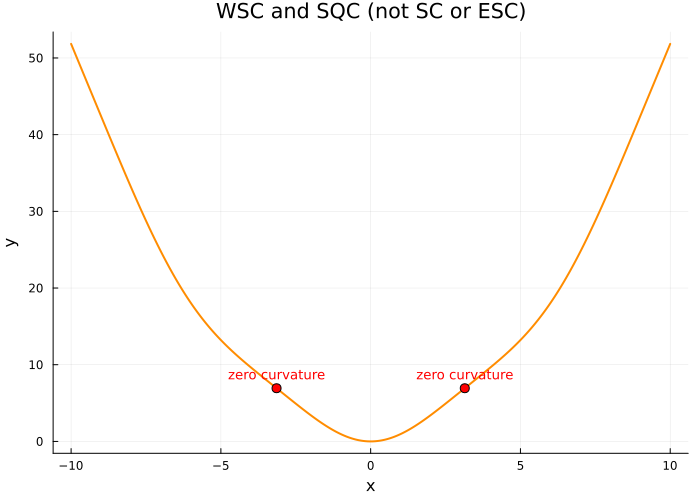}
\includegraphics[width=.4\textwidth]{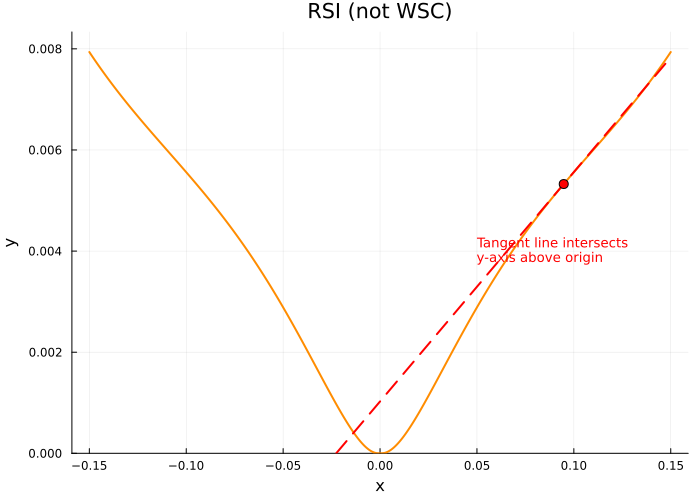}
\includegraphics[width=.45\textwidth]{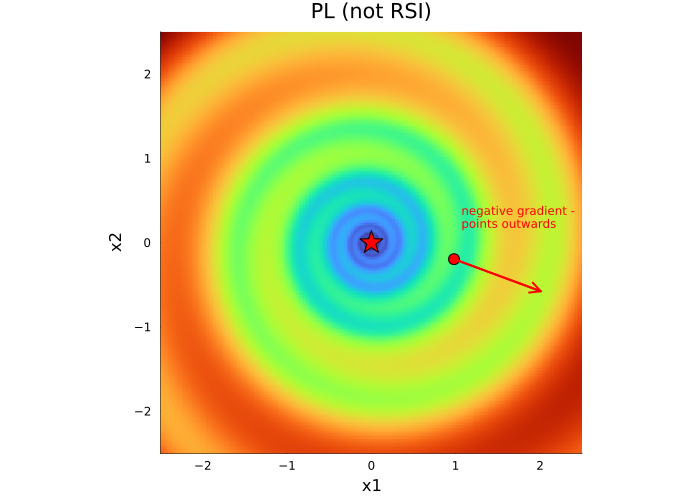}
\includegraphics[width=.45\textwidth]{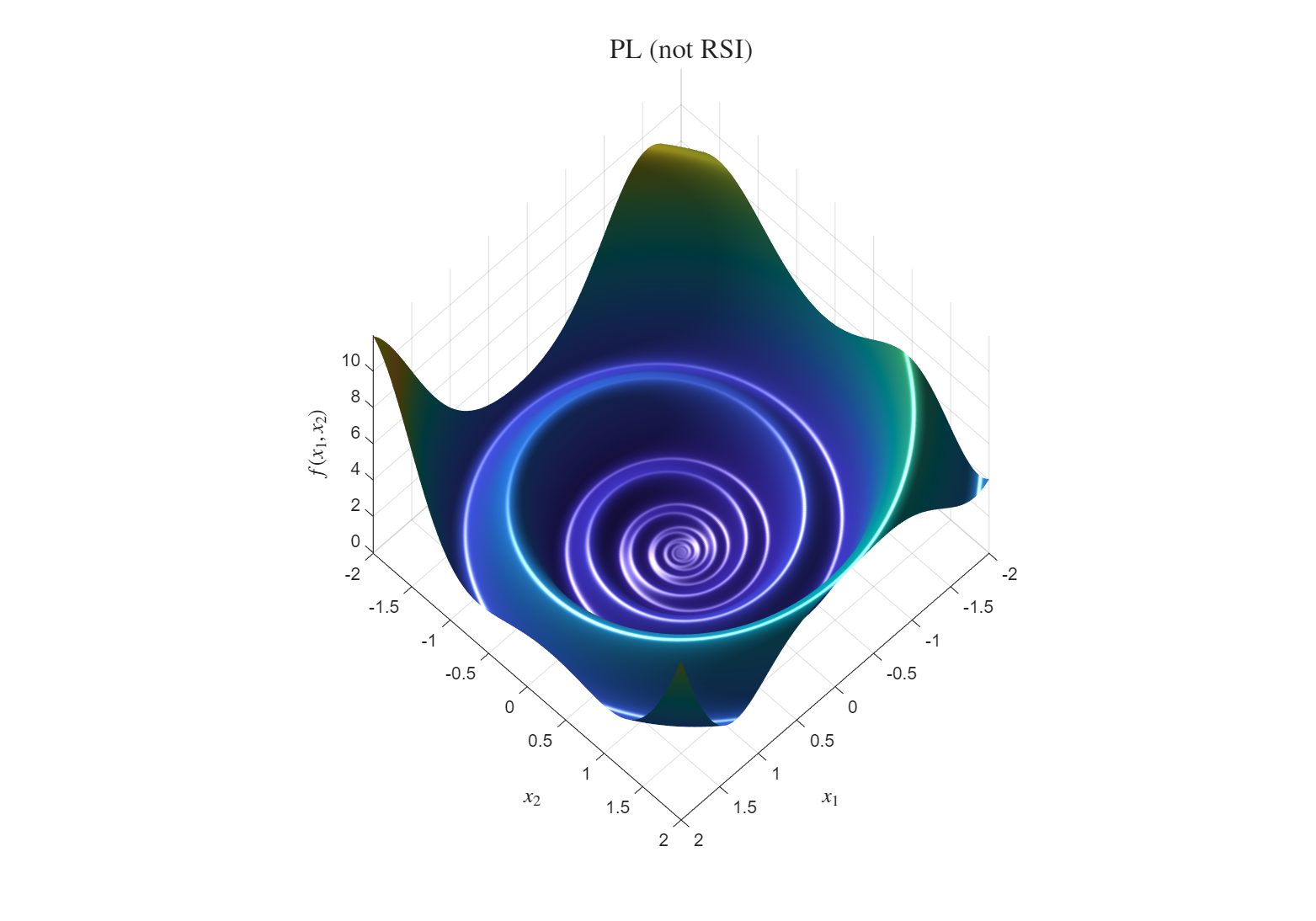}
\includegraphics[width=.4\textwidth]{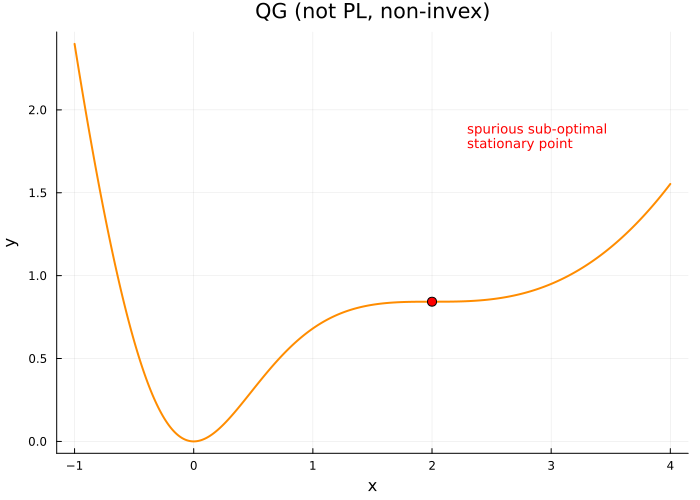}
\includegraphics[width=.4\textwidth]{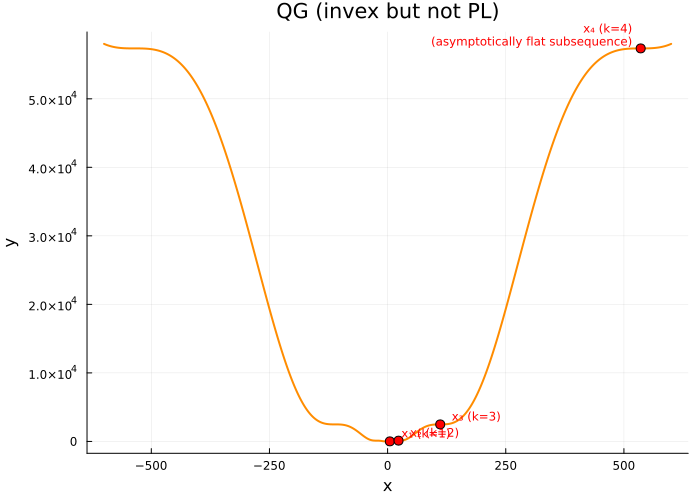}
\caption{A strongly-convex function (upper left), and counter-examples satisfying one condition but not another.}
\label{fig:counter}
\end{center}
\end{figure}

Below we give the updated Gemini-assisted result relating the conditions for smooth functions. As before, all of these conditions except QG are sufficient to show linear convergence of gradient descent.
\begin{theorem}\label{thm:2}
For a function $f$ with a Lipschitz-continuous gradient, the following implications hold:
\[
\begin{array}{rcccl}
& & (ESC) \subset (WSC) & & \\
& \ru{\subset} & & \rd{\subset} & \\
(SC) & & & & (QSC) \equiv (RSI) \subset (EB) \equiv (PL) \subset (QG) \\
& \rd{\subset} & & \ru{\equiv} & \\
& & (SQC) \subset (SSQC) & &
\end{array}
\]
These  relationships continue to hold if we assume that $f$ is invex. We further obtain the following equivalences under additional assumptions:
\begin{itemize}
\item If $f$ is convex, then we have 
\[
\begin{array}{rcccl}
& & (ESC) \subset (WSC) & & \\
& \ru{\subset} & & \rd{\equiv} & \\
(SC) & & & & (QSC) \equiv (RSI) \equiv (EB) \equiv (PL) \equiv (QG) \\
& \rd{\subset} & & \ru{\equiv} & \\
& & (SQC) \subset (SSQC) & &
\end{array}
\]
\item If $f$ is convex and quadratic, then we have 
\[
\begin{array}{rcccl}
& & (ESC) \equiv (WSC) & & \\
& \ru{\subset} & & \rd{\equiv} & \\
(SC) & & & & (QSC) \equiv (RSI) \equiv (EB) \equiv (PL) \equiv (QG) \\
& \rd{\equiv} & & \ru{\equiv} & \\
& & (SQC) \subset (SSQC) & &
\end{array}
\]
\end{itemize}
\end{theorem}
\vspace{-1em}
We give the counter-examples and additional implications in Appendix~\ref{sec:A}. Appendix~\ref{sec:A} also shows that SC and ESC are equivalent if there is a unique solution, and that RSI and PL are equivalent for univariate functions. Figure~\ref{fig:counter} plots a selection of the counter-examples. Note that there was only one case where the counter-example was required to have more than one dimension: PL implies RSI in one dimension so we required a two-dimensional ``toilet flush'' function as a counter-example satisfying PL but not RSI. We found the ``toilet flush'' particularly interesting: the negative gradient can point in a direction away from the optimal solution yet gradient descent still converges linearly to the solution due to the PL property.

We note that~\citet{Jovanovic_1993} gives the equivalence between SC and SQC for convex quadratic functions, and that~\citet{khanh2025star} and~\citet{khademi2026convexity} had already shown the connections between SSQC, SQC, and RSI. 
We also note that~\citet{hui} present a recent variant of this hierarchy with two additional conditions (RSC' and PL'), while~\citet{khademi2026convexity} gives relationships between various conditions related to convexity and invexity (not necessarily related to linear convergence of GD).
Further equivalences are explored for $C^2$ functions near local minima by~\citet{rebjock2025fast} who show that QG is equivalent to EB and PL in this setting, as well as the Morse-Bott property. We note that to construct an $L$-smooth function satisfying QG but not PL, we either needed to violate invexity or needed to create a function that violates PL in the limit as we move away from the minimizers. For the issue of whether QG implies PL for smooth invex functions over a compact set containing the solutions, we used Gemini to obtain incorrect proofs in both directions.

\section{Criticisms}
\label{sec:crit}

There have been a variety of criticisms of our paper and the PL condition itself, dating back to even before the paper was published. We summarize some of the main criticisms below:
\begin{itemize}
\item {\bf It can be hard to verify that a function satisfies PL.} Unlike strong-convexity, we are not aware of a simple way to characterize that a function satisfies PL. For example, for strong-convexity we may be able to lower bound the Hessian matrix globally, or we can use that the sum of a strongly-convex function and a convex function is strongly-convex. However, the sum of two PL functions is not necessarily PL. Indeed, when we first discussed the work with Matt Hoffman (the DeepMind instantiation) at ICML 2015, he asked if this was a variation on using the~\citet{hoffman} (no relation) constant. While the PL condition theoretically applies to a wide variety of possible functions, to show that it holds for least squares and many other problem classes, we do indeed use the Hoffman constant. \smallskip
\item {\bf The PL constant $\mu$ can be arbitrarily small.} When we presented the work at the 2015 CIFAR workshop, Francis Bach highlighted that the PL constant may be so small that the linear convergence rate may be slower than the sublinear rate of gradient descent for convex functions (for any reasonable finite number of iterations). In the case of least squares the PL constant $\mu$ is the smallest non-zero singular value of the Hessian, which can be arbitrarily close to zero. \smallskip
\item {\bf PL is equivalent to EB.} Ruoyu Sun (who was a student of Tom Luo, one of the authors of the EB paper) discussed the work with us at the 2015 NeurIPS Optimization for ML workshop. He did not seem to think the paper was interesting, as the PL condition seemed to hold for a similar class of functions as the EB condition. Ruoyu ended up being correct, as in 2016 it was established that EB and PL are equivalent, even in the general non-convex setting. \smallskip
\item {\bf Missed/mis-cited results, particularly the work of Jérôme Bolte.} The ECML version missed and/or mis-cited several results. For example, while we cited~\citet{Zhang2015} as prior work on the relationships between conditions in the convex setting, it was missing a citation to~\citet{Bolte2015} for the same result. It also did not include a citation to the concurrent~\citet{zhang2016characterization}, who also considered the general setting. The ECML version also incorrectly stated that the KL condition had only previously been used to show local convergence of proximal-gradient methods, whereas works like~\citet{Bolte2015} had indeed shown that the KL condition leads to a global linear rate. The ECML version was also missing the equivalence between the proximal-PL and KL conditions. Regrettably, both Hui Zhang and Jérôme Bolte had to email us to bring these oversights to light. We promptly corrected the arXiv version of the paper in each instance (v2 for Hui Zhang, v3 for Jérôme Bolte). \smallskip
\item {\bf Incorrect results}. The original arXiv paper (v1) incorrectly states that the proximal-PL condition holds if $f$ has the form $f(x)=h(Ax)$ for a strongly-convex $h$ and if $g$ is convex (not necessarily polyhedral). However, without the polyhedral assumption, the proximal-gradient method does not converge linearly in general. Several people contacted us about this (Ting Kei Pong and Zirui Zhou), and the AI/Stats reviewers also caught this error, which was a contributing factor to the paper being rejected. This was fixed in the arXiv paper (v2) and in the ECML version. We later updated the paper again (v4) in 2020 when Stephen Becker and Liam Madden pointed out a missed (1/2) factor in the relationships between conditions that ended up simplifying the result. \smallskip
\item {\bf Stochastic gradient descent (SGD) assumptions are contradictory.} To control the variation in the stochastic gradient estimates, the stochastic gradient descent results in the paper assume $\mathbb{E}[\norm{\nabla f_i(x)}^2] \leq \sigma^2$ for some $\sigma^2 > 0$. However, globally this assumption contradicts the PL condition which requires the gradient to grow without bound. These assumptions can hold if the iterates are constrained to a bounded set, or alternately the results can be modified to use the assumption that $\mathbb{E}[\norm{\nabla f_i(x) - \nabla f(x)}^2] \leq \sigma^2$. We have used the latter result in class but we have not seen it in the literature except as special cases of more-general/complicated results, so we include it in Appendix~\ref{sec:B}. We also note that our variance-reduced SGD assumed convexity, and works not requiring convexity for variance-reduced SGD appeared shortly after~\citep{reddi2,reddi1}.
 \smallskip
\item {\bf PL functions are ``basically'' SC}. This is true in the sense that PL functions act like SC functions if we follow gradient trajectories away from the solution. Thus, even though they can technically be non-convex, they are quite constrained in how the non-convexity can appear. Related to this are the recent results  that $C^\infty$ globally PL functions are equivalent to nonlinear least-squares problems (linear least squares up to a diffeomorphism), and $C^2$ PL functions with bounded argmins are locally SC around the solution~\citep{nejma2025polyak,boumal2026smooth}. \smallskip
\item {\bf PL does not allow Nesterov acceleration.} While PL can replace SC in a wide variety of settings,~\citet{yue2023lower} show that PL is insufficient to achieve the rate of Nesterov acceleration on SC functions. We discuss this issue further in the next section. \smallskip
\item {\bf PL cannot hold globally for neural networks.} The PL condition cannot hold for standard neural network training objective functions because they are not invex. We discuss this issue further in the next section.
\end{itemize}

\section{Related and Follow-Up Works}
\label{sec:ext}

There have been a wide variety of follow-up works in the 10 years since the paper was first released. Here we summarize a variety of notable follow-ups, and also highlight several related works that appear before our paper.
\begin{itemize}
\item {\bf Iterate convergence}: In our work we focus on the convergence of the function values $f(x_k)$ where $x_k$ are the iterations of GD with a suitable step size. But one might also be curious about convergence of the gradient descent iterates $x_k$ themselves. For smooth PL functions it follows from the equivalence of PL and EB that  $\norm{x_k - (x_k)_p}^2$ converges to zero when we run GD with suitable step sizes. But the original work of~\citet{polyak} shows the stronger property that the sequence $\norm{x_k - x_*}^2$ converges with a linear rate to a particular solution $x_*$ (not necessarily the projection of $x_0$ onto the solution).
\item {\bf Weak PL condition}: This is a condition introduced by~\citet{csiba2017global} allowing non-convexity but preserving the convergence rate of gradient descent on convex functions rather than strongly-convex functions. This is a generalization of the earlier star-convex property defined by~\citet{nesterov2006cubic} and of the later quasar convexity~\citep{hinder2020near}. In particular, the hierarchy between these conditions is that star convexity implies quasar convexity which implies weak PL (analogous to how WSC implies QSC which implies PL). \smallskip
\item {\bf Gradient dominated functions}: This is a generalization of the PL condition where the exponent on the gradient norm can take a value other than 2. It was defined by~\citet{nesterov2006cubic} who used it (and star convexity) to analyze Newton's method with cubic regularization. \citet{wilson2019accelerating} show that gradient descent with adaptive step sizes can achieve a linear rate for gradient-dominated functions with exponent greater than 2. \smallskip
\item {\bf Kurdyka-\L{}ojasiewicz condition}: This is a generalization of gradient-dominated functions (and hence PL) that allows the function to be non-smooth (beyond the proximal-gradient setting we consider). It was developed by~\citet{bolte2010characterizations} building on the work of~\citet{kurdyka1998gradients}. The KL condition leads to convergence rates for a more general class of functions than the PL and proximal-PL conditions. Indeed,~\citet{attouch2013convergence} analyze the convergence of a variety of algorithms including for non-smooth and non-convex problems (possibly with inexact computation) and even  non-continuous problems like $\ell_0$ minimization.
The EB condition can also be used in non-smooth settings, and for convex functions the KL and EB conditions are equivalent~\citep{Bolte2015,drusvyatskiy2016error}. \smallskip
\item {\bf Tight analyses under PL}: \citet{abbaszadehpeivasti2023conditions} use the performance estimation problem (PEP) framework to analyze GD under PL. They derive a faster linear convergence rate and optimal step sizes. They also argue in their Theorem~5 that the PL assumption is a necessary assumption for linear convergence of GD.\smallskip
\item {\bf Stochastic gradient descent under interpolation}: If the noise converges to zero as we approach the solution, stochastic gradient descent can obtain a linear convergence rate on strongly-convex functions similar to gradient descent. \citet{vaswani2019fast} show that this property is preserved under the PL condition. \smallskip
\item {\bf Stochastic gradient descent with random reshuffling}: Several works have used the PL condition to analyze the variant on stochastic gradient descent where, instead of sampling independent random examples on each iteration, we go through random permutations of the examples~\citep{haochen2019random,ahn2020sgd,mishchenko2020random,nguyen2021unified}. This variant has also been studied under interpolation assumptions~\citep{fan2023fast}. \smallskip
\item {\bf Accelerated gradient methods}: While~\citet{yue2023lower} argue that accelerated convergence rates via Nesterov-style algorithms are not possible under only the PL condition, several authors have considered additional assumptions that allow acceleration. For example,~\cite{hermant2026acceleration} proposes an additional aiming condition that allows acceleration, while~\citet{feng2026optimal} show that local acceleration near the solution is possible for $C^2$ functions. \smallskip
\item {\bf Policy gradient methods}: Several works consider using PL-style inequalities in reinforcement learning settings to analyze policy gradient methods~\citep{fazel2018global,mei2021leveraging}. \smallskip
\item {\bf Riemannian PL}: A variety of works consider local PL conditions defined on manifolds. Many classic problems such as principal component analysis can be analyzed in this framework. A recent work discussing this topic is~\citet{wang2025distributed}. \smallskip
\item {\bf Log-Sobolev inequality}: This condition can be viewed as an analogue of PL for probability spaces~\citep{diaconis1996logarithmic}. It has been used to analyze Langevin sampling~\citep{vempala2019rapid}, expectation maximization~\citep{caprio2025fast}, and sampling in diffusion models~\cite{lee2023convergence}. \citet{chewi2024ballistic} discuss how the zero-temperature limit of the log-Sobolev constant is the PL constant, while \citet{blanchet2018family} discuss how the log-Sobolev inequality is a \L{}ojasiewicz inequality discussing connections to several classical functional inequalities. \smallskip
\item {\bf Local PL condition for optimal transport}: \citet{gonzalez2026polyak} argue that the dual problem associated with quadratically-regularized optimal transport satisfies a PL condition locally around the solution. \smallskip
\item {\bf Local PL conditions on neural networks}: While neural network training objectives tend to be non-invex and thus do not satisfy the PL condition, a variety of works have shown that local PL-like conditions can hold for neural networks. The early works on this topic required strong assumptions like the ability to interpolate the data, large widths, and the use of the squared error~\citep{allen2019convergence}. However, through a sequence of works, these assumptions have all been relaxed to argue that local PL variants may help explain the success of gradient methods in training neural networks~\citep{liu2020toward,frei2021proxy,barboni2022global,scaman2022convergence,aich2025sublinear}. \smallskip
\item {\bf PLNets}: \citep{wang2024monotone} introduce PLNets, which are neural networks designed to satisfy the PL condition with respect to their inputs allowing efficient optimization of inputs. \smallskip
\end{itemize}
Finally, the third author would like to shamelessly plug a related recent work out of their lab. This work started with the question:
\begin{quote}``\emph{What is the simplest assumption under which gradient descent with line search can converge arbitrarily faster than with fixed step sizes?}''\end{quote}
This led to our work on the glocal smoothness assumption, which is written in a similar spirit to the PL paper~\citep{fox2025glocal}. 
While the peer review process tends to reward increasingly-complicated mathematization, we encourage the reader to also consider questions about \emph{simplicity} that could benefit the literature.

\section*{Acknowledgments} 

We are grateful to several people for providing quick feedback on this hastily-written report. In particular we would like to thank J\'{e}rome Bolte, Nicolas Boumal, Robert Gower, Abbas Khademi, Felipe Lara, Gabriel Peyr\'{e}, and Hui Zhang. This work was supported by the Canada CIFAR AI Chairs Program and NSERC Discovery Grant RGPIN 2022-036669.

\section*{References}
\vspace{-2em}
\renewcommand{\bibname}{}
\bibliographystyle{abbrvnat}
\bibliography{bibGemini}

\appendix

\section{Only SC Implies Convexity}
\label{sec:convex}

In this section, we clarify several relationships between the conditions and convexity for differentiable functions (not necessarily smooth), assisted by Gemini. The below implications yield that the only condition implying convexity is SC, except in 1 dimension where ESC also implies convexity.

\textbf{ESC implies convexity for univariate functions}. Since ESC implies SC if the solution is unique, it remains to examine the case where we have an interval $[a,b]$ of solutions. We can divide the space into 3 regions: the region $(-\infty,a]$ projects onto $a$ and thus ESC implies SC on this region, the region $[b,\infty)$ projects onto $b$ and thus ESC implies SC on this region, while the interval $[a,b]$ is constant and hence convex. Between regions, the convexity inequality is satisfied since the derivative is negative on $(-\infty,a]$, 0 on $[a,b]$, and positive on $[b,\infty)$.

\textbf{ESC does not imply convexity}. Let $f(x_1, x_2) = x_1^2(2 + \sin(x_2))$. The solution set is $x_p=(0,x_2)$ for any $x_2$ so if $x$ and $y$ share the same projection they have the same $x_2$ coordinate and the function behaves like $f(x_1) = x_1^2(2+\sin(x_2))$ which is SC implying ESC. The second derivative with respect to $x_2$ is $-x_1^2\sin(x_2)$ which is negative at $x_1=1$ and $x_2=\pi/2$ implying the function is non-convex.

\textbf{WSC does not imply convexity for univariate functions}. Let $f(x)=x^4 +x^3 + 0.35x^2$. Its second derivative  $f''(x)=12x^2 + 6x + 0.7$ is $-1/20$ at $x=-1/4$ hence the function is non-convex. The function achieves its minimum of 0 at $x=0$. The function $f$ satisfies WSC with $\mu=0.02$ since $f(x) + f'(x)(0-x) + 0.01x^2 = (x^4 + x^3 + 0.35x^2) - x(4x^3 + 3x^2 + 0.7x) + 0.01x^2 = -x^2(3x^2 + 2x + 0.34)$ and the right factor is positive (so the expression is at most 0).

\textbf{SQC does not imply convexity for univariate functions}. Let $f(x) = x^6 - 3x^4 + 4x^2$. Its second derivative  $f''(x) = 30x^4 - 36x^2 + 8$ is $-14/5$ at $x=\sqrt{3/5}$ hence the function is non-convex. To show that it satisfies SQC we need to show that $f'(y)(y-x) \geq \frac \mu 2(y-x)^2$ for some positive $\mu$ whenever $f(y)\geq f(x)$. Note that $f$ is strictly increasing in $|x|$ so $f(y) \geq f(x)$ is equivalent to $|y|\geq |x|$. If $y \geq 0$ then $f(y) \geq f(x)$ implies $x\in[-y,y]$  and thus $(y-x)\leq 2y$. We also have $f'(y) = 6y^5 - 12y^3 + 8y = 2y(3y^4 - 6y^2 + 4)\geq 2y$ since $3y^4 -6y^2 + 4 = 3(y^2 - 1)^2 + 1 \geq 1$. Using $f'(y)\geq 2y$ and $2y \geq (y-x)$ we have $f'(y)(y-x)\geq (y-x)^2$ so SQC holds with $\mu=2$ for $y>0$. By the symmetry of the function around 0 it also holds for $y<0$.

\section{Relationships Between Conditions (Again)}
\label{sec:A}

In this section, we give a variety of implications and counter-examples regarding the conditions, obtained through interacting with Gemini. Together with the previous 2016 result in Section~\ref{sec:thm2}, these imply the updated result in Section~\ref{sec:updated}. We assume differentiability throughout this section, but highlight cases where we additionally require smoothness.

\textbf{SC implies SQC}. Re-arranging SC we have $0 \geq f(x) - f(y) + \langle \nabla f(x),y -x\rangle +\frac \mu 2\norm{y-x}^2$. If $f(y)\leq f(x)$ then $f(x) - f(y)\geq 0$ which yields SQC with the same $\mu$.

\textbf{ESC implies SC if the solution is unique}. If there is a single global minimum $x^*$, then for every $x$ and $y$ we have $x_p=y_p$ so the ESC condition implies the SC condition for all $x$ and $y$.

\textbf{ESC does not imply SC}. Let $f(x) = \half\text{dist}(x,[-1,1])^2$ which has an interval of solutions, $[-1,1]$. Since the function has $f''(x)=0$ for $x$ in $[-1,1]$ it does not satisfy SC. 

ESC requires SC between points with the same projection onto the solution set. The interval $(-\infty,-1]$ projects onto $-1$ and is SC while the interval $[+1,\infty)$ projects onto $+1$ and is SC. If $x$ and $y$ are in $[-1,1]$ then both $x$ and $y$ are solutions so to have the same projection we must have $x=y$ and the ESC condition is satisfied ($0 \geq 0$).

Note that this example also does not satisfy SQC since the solution is not unique, so this example also shows that ESC does not imply SQC.

\textbf{ESC does not imply SC for convex quadratic functions}. Let $f(x_1,x_2)=\half x_1^2$. This function is flat along $x_2$ so does not satisfy SC. The solution set is $x_p=(0,x_2)$ for any $x_2$ so if $x$ and $y$ share the same projection they have the same $x_2$ coordinate and the function behaves like $f(x_1) = \half x_1^2$ which is SC.

\textbf{WSC implies QSC}. If the WSC inequality holds then this is exactly the QSC inequality with $\gamma=1$.

\textbf{WSC does not imply ESC for convex functions}. Let $f(x) = \half x^2 - \cos(x) +1$ which has a unique optimum of $0$ at $x=0$ and $f''(x) = 1+\cos(x)\geq 0$. We have $f''(\pi)=0$ so the ESC inequality cannot be satisfied at this point. 

The WSC inequality requires $0 \geq (\half x^2 - \cos(x) + 1) - x(x +  \sin(x)) + \frac \mu 2 x^2$. Using $\mu=0.2$ in the right side gives $-0.4x^2 - x\sin(x) - \cos(x)  + 1$ which obtains a maximum value of 0 at 0.

\textbf{WSC implies ESC for convex quadratic functions}. Let $f(x) = \half x^TAx + b^Tx +c$ with $A \succeq 0$. Using that $\nabla^2f(x)=A$ we get from a Taylor expansion that $f^* = f(x) + \langle \nabla f(x),x_p - x\rangle + \half (x-x_p)^TA(x-x_p)$. Using this in WSC gives for all $x$ that $\half (x-x_p)^TA(x-x_p) \geq \frac \mu 2\norm{x_p-x}^2$. Consider an $x$ and $y$ that share the same projection $x_p$. To obtain the ESC inequality define $w = x_p + (y-x)$. Note that we have $w_p = x_p$ since  $y-x = (y-x_p) - (x-x_p)$ and both terms are orthogonal to the null space of $A$. Evaluating WSC at $w$ thus gives $\half (y-x)^TA(y-x) \geq \frac \mu 2\norm{y-x}^2$ and we obtain ESC by using this to bound the Taylor expansion: $f(y) = f(x) + \langle \nabla f(x),y -x\rangle + \half(y-x)^TA(y-x) \geq f(x) + \langle \nabla f(x),y- x\rangle + \frac \mu 2 \norm{y-x}^2$. Combined with ESC implying WSC in general, this establishes the equivalence of ESC and WSC for convex quadratic functions.

\textbf{SQC implies SC for convex quadratic functions}. Let $f(x) = \half x^TAx + b^Tx +c$ with $A \succeq 0$. For the purpose of contradiction assume there exists a $v \neq 0$ such that $Av=0$. Note that this implies that $\langle \nabla f(x),v\rangle = \langle Ax+b,v\rangle = b^Tv$. If $b^Tv=0$ then $f(x+v)=f(x)$ and SQC requires $0 \geq \langle Ax + b,v\rangle + \frac \mu 2 \norm{v}^2 = 0 + \frac \mu 2 \norm{v}^2$ which is a contradiction since $ v\neq 0$. If $b^Tv > 0$ consider any $x$ and $y$ such that $y=x-\alpha v$ so that $f(y) = f(x) - \alpha b^Tv < f(x)$. SQC requires $0 \geq -\alpha\langle \nabla f(x),v\rangle + \frac{\mu\alpha^2}{2}\norm{v}^2 = -\alpha b^Tv + \frac{\mu\alpha^2}{2}\norm{v}^2$ which cannot be satisfied for sufficiently large $\alpha$. A similar argument holds for $b^Tv < 0$ and thus we have that $A$ is positive definite and the function is SC. Combined with SC implying SQC in general, this establishes the equivalence of SC and SQC for convex quadratic functions.

\textbf{SQC does not imply SC for convex functions}. Let $f(x) = \half x^2 - \cos(x) +1$ which we showed above is convex but not ESC and thus not SC.

To show that it satisfies SQC we need to show that $f'(y)(y-x) \geq \frac \mu 2(y-x)^2$ for some positive $\mu$ whenever $f(y)\geq f(x)$. Note that $f$ is strictly increasing in $|x|$ so $f(y) \geq f(x)$ is equivalent to $|y|\geq |x|$. If $y \geq 0$ then $f(y) \geq f(x)$ implies $x\in[-y,y]$  and thus $(y-x)\leq 2y$. We have $f'(y) = y + \sin(y) = y\left(1 + \frac{\sin(y)}{y}\right) \geq 0.5y$ for $y \geq 0$. Using $f'(y) \geq 0.5y$ and $2y \geq (y-x)$ we have $f'(y)(y-x)\geq 0.25(y-x)^2$ so SQC holds with $\mu=0.5$ for $y>0$. By the symmetry of the function around 0 it also holds for $y<0$.

Note that since this example is also not ESC, it also shows that SQC does not imply ESC.

\textbf{SQC implies SSQC}. For any $x$, set $y=x_p$ in the SQC inequality and use that $\max\{f(x),f(x_p)\} = f(x)$ to obtain the SSQC inequality.

\textbf{SSQC implies RSI}. Defining $\alpha=1-\lambda$ and re-arranging SSQC gives $\frac{f(x) - f(x + \alpha(x_p - x))}{\alpha} \geq \frac{\mu}{2}(1-\alpha)\norm{x_p - x}^2$. Taking the limit as $\alpha\to 0^+$ and using that the left side limit is $\langle \nabla f(x), x - x_p\rangle$ gives RSI with a constant of $\mu/2$.

\textbf{QSC implies RSI}. Re-arranging QSC gives $\frac 1 \gamma \langle \nabla f(x),x-x_p\rangle \geq f(x) - f^*  + \frac \mu 2 \norm{x_p-x}^2$. Using that $f(x)-f^*\geq 0$ and multiplying by $\gamma$ gives RSI with a constant of $\gamma\mu/2$. 

\textbf{QSC does not imply WSC}. Let $f(x) = x^2(1.3 + \sin(\log|x|))$ but define $f(0)=0$. This has a unique minimum of $0$ at $x=0$, so WSC requires $0 \geq f(x) - xf'(x) + \frac \mu 2 x^2$ or equivalently $\frac \mu 2 \leq \frac{xf'(x) - f(x)}{x^2}$. We have $f'(x) = 2.6x + 2x\sin\log|x| + x\cos\log|x|$ and thus the WSC requirement simplifies to $\frac \mu 2 \leq 1.3 + \sin\log|x| + \cos\log|x|$. Choosing $x$ such that $\log|x| = -\frac{3\pi}{4}$ gives $\frac \mu 2 \leq 1.3 - \sqrt{2} \approx -0.1142$ which cannot hold for any positive $\mu$.

For QSC we only require $\frac \mu 2 \leq \frac{\frac 1 \gamma xf'(x) - f(x)}{x^2}$ which is equivalent to $\frac \mu 2 \leq \frac 1 \gamma(2.6 + 2 \sin\log|x| + \cos\log|x|) - (1.3 + \sin\log|x|)$ = $\left(\frac{2.6}{\gamma} - 1.3\right) + \left(\frac 2 \gamma - 1\right)\sin\log|x| + \frac 1 \gamma \cos\log|x|$. With $\gamma=0.3$ the right side is minimized at approximately $0.792$ (using the harmonic addition identity that $\theta_1\sin(x) + \theta_2 \cos(x) \geq -\sqrt{\theta_1^2+\theta_2^2}$) and we can choose $\mu =1.5$ giving $\mu/2=0.75$ to satisfy QSC.

\textbf{RSI implies QSC for smooth functions}. Re-arranging RSI and multiplying by $1/\gamma$ gives
\[
\frac 1 \gamma \langle \nabla f(x), x - x_p\rangle \geq \frac \mu \gamma \norm{x - x_p}^2.
\]
Adding this to the descent lemma from L-smoothness,
\[
f(x) \leq f^* + \frac L 2\norm{x-x_p}^2
\]
gives after re-arranging that
\[
f^* \geq f(x) + \frac 1 \gamma \langle \nabla f(x), x_p - x\rangle + \left(\frac \mu \gamma - \frac L 2\right)\norm{x-x_p}^2.
\]
Choose $\gamma = \mu/L$ to give
\[
f^* \geq f(x) + \frac L \mu \langle \nabla f(x), x_p - x\rangle + \frac L 2\norm{x-x_p}^2,
\]
which is QSC with $\gamma = \mu/L$ and $\mu=L$ (note that QSC allows $\mu$ to be larger than $L$ while many conditions like RSI require their $\mu$ to be at most $L$). This establishes the equivalence of RSI and QSC when combined with QSC implying RSI.

\textbf{RSI implies SSQC}. Define $x_\alpha = \alpha x + (1-\alpha)x_p$ so that $x_\alpha - x_p = \alpha(x - x_p)$ and apply RSI to get
\[
\langle \nabla f(x_\alpha), x-x_p\rangle \geq \mu \alpha\norm{x - x_p}^2.
\]
Define the function evaluated along the path between $x$ and $x_p$ as $g(\alpha) = f(
\alpha x + (1-\alpha)x_p)$ so $g'(\alpha) = \langle \nabla f(x_\alpha), x - x_p\rangle$. We then have for $\lambda \in [0,1]$ that
\[
f(x) - f(x_\lambda) = g(1) - g(\lambda) = \int_\lambda^1 g'(\alpha)d\alpha \geq \int_\lambda^1\mu\alpha\norm{x-x_p}^2d\alpha = \frac \mu 2 (1-\lambda^2)\norm{x-x_p}^2.
\]
The SSQC inequality follows from re-arranging and using that $1-\lambda^2 \geq \lambda(1-\lambda)$ on $[0,1]$. This establishes the equivalence of SSQC and RSI when combined with SSQC implying RSI, and further establishes the equivalence of SSQC with QSC due to the equivalence of QSC and RSI.

\textbf{RSI does not imply SQC for convex quadratic functions}. Both of our above ``ESC does not imply SC'' examples also show that RSI does not imply SQC: they satisfy ESC (hence also RSI) and have multiple solutions (hence do not satisfy SQC). Both of these functions are convex and one is quadratic. Note that this establishes the separation of SQC and SSQC due to the equivalence of RSI and SSQC.

\textbf{EB implies WSC for smooth convex functions}: The co-coercivity property of convex $L$-smooth functions is that 
\[
f(y) \geq f(x) + \langle \nabla f(x),y-x\rangle + \frac{1}{2L}\norm{\nabla f(x)-\nabla f(y)}^2.
\]
Choosing $y=x_p$ and using EB gives
\[
f(x_p) \geq f(x) + \langle \nabla f(x),x_p-x) + \frac{\mu^2}{2L}\norm{x_p-x}^2,
\]
which is WSC with constant $\mu^2/L$. Combined with WSC implying EB (in the original paper), this establishes the equivalence of WSC and EB for smooth convex functions.

\textbf{PL does not imply RSI}. Consider the function $f(r,\theta) = r^2(1 + 0.5\sin(\theta - 10\log r))$ where $x \in \R^2$ is represented in polar coordinates. This has a unique minimum of $0$ at $x=(0,0)$. The left-side of the RSI inequality is thus $\langle \nabla f(x), x - x_p\rangle = \langle \nabla f(x), x \rangle = r\frac{\partial f}{\partial r} = r^2(2 + \sin \alpha - 5 \cos \alpha)$ for $\alpha = \theta - 10\log r$. 
To show that RSI is not satisfied choose $r=1$ so that $\log r =0$ and $\theta = \arctan(-0.2)$ so that $\sin\alpha=\sin\theta=-1/\sqrt{26}$ and $\cos\alpha = \cos\theta=5/\sqrt{26}$, corresponding to $x \approx (0.9806, -0.1961)$. This point gives a left-side of $1^2\left(2 - \frac{1}{\sqrt{26}} - 5\frac{5}{\sqrt{26}}\right) = 2 -\sqrt{26} \approx -3.099$. Since this is negative (the gradient points away from the solution), no positive $\mu$ can make the RSI inequality satisfied.

Showing that the function satisfies PL is equivalent to showing that $\frac{\norm{\nabla f(x)}^2}{f(x)}\geq 2\mu$ for some positive $\mu$. The left side simplifies to $\frac{(2+\sin\alpha - 5\cos\alpha)^2 + 0.25\cos^2\alpha}{1+0.5\sin\alpha}$. The denominator is between 0.5 and 1.5, while the two non-negative terms in the numerator cannot simultaneously approach 0: to make the second zero we require $\cos\alpha$ to approach 0 which means that $\sin\alpha$ approaches $1$ or $-1$ making the first term converge to 1 or 9.

\textbf{PL implies RSI for univariate functions}. Take some $x > x_p$ so that $f'(x) > 0$. Recall that PL implies QG with the same $\mu$. Using PL and QG we have $f'(x) \geq \sqrt{2\mu(f(x)-f^*)} \geq \sqrt{2\mu\left(\frac \mu 2(x-x_p)^2\right)} = \mu(x-x_p)$. Multiplying by $(x-x_p)$ gives RSI with the same $\mu$, and repeating the steps for $x < x_p$ gives the result.

\textbf{QG does not imply PL}. Let $f(x) = \half x^2 - 4x + \frac 3 2 \log(x^2+1) + 4 \arctan(x)$. This has a unique global minimum of $0$ at $x=0$. The derivative is $f'(x) = \frac{x(x-2)^2}{x^2+1}$. We have $f'(2)=0$ so the function has a stationary point that is not a global minimum and hence cannot satisfy PL.

Showing that the function satisfies QG is equivalent to showing that $\frac{f(x)}{x^2} \geq \frac \mu 2$ for some positive $\mu$. The function $\frac{f(x)}{x^2}$ is minimized at $x\approx 3.65$ with a value of $0.0951$ so we can take $\mu=0.01$.

\textbf{QG does not imply PL for invex functions}. Let $f(x) = \frac{\exp(u(x))}{20}(5 - \cos(2u(x)) - 2\sin(2u(x))) + \frac{u(x)}{2} - \frac 1 5$ with $u(x) = \log(x^2+1)$. This has a unique global minimum of $0$ at $x=0$. The derivative simplifies to $f'(x) = x\sin^2(u(x)) + \frac{x}{x^2+1}$. Consider the sequence defined by $x_k = \sqrt{\exp(k\pi)-1}$: we have $u(x_k) = k\pi$ and $\sin(k\pi)=0$ so we have
$f(x_k) = \frac{\exp(k\pi)}{20}(5 - 1 - 0) + \frac{k\pi}{2} - \frac 1 5 = \frac 1 5(\exp(k\pi)-1) + \frac{k\pi}{2}$ and $f'(x_k) = \frac{x_k}{x_k^2+1} \leq \frac{\sqrt{\exp(k\pi)}}{\exp(k\pi)} = \frac{1}{\sqrt{\exp(k\pi)}}$. We thus have $\frac{(f'(x_k))^2}{2f(x_k)} \leq \frac{1}{2\exp(k\pi)(\frac 1 5(\exp(k\pi)-1) + \frac{k\pi}{2})}$ which converges to zero as $k$ grows and thus no positive $\mu$ bounds the ratio.

Showing that $f$ satisfies QG is equivalent to showing $\frac{f(x)}{x^2} \geq \frac \mu 2$ for some positive $\mu$ for $x\neq 0$. Using $\cos(2u)+2\sin(2u)\leq \sqrt{1^2 + 2^2} = \sqrt{5}$ we have $5-\cos(2u(x))-2\sin(2u(x)) \geq 5 -\sqrt{5}$. Using that $u(x)\geq 0$ we have $f(x) \geq \frac{5-\sqrt{5}}{20}(x^2+1) - \frac 1 5$. For $|x|\geq 2$ we thus have $\frac{f(x)}{x^2} \geq \frac{5-\sqrt{5}}{20} + \frac{1-\sqrt{5}}{80} > 0.1$. For $x \in [-2,2]$, defining $\frac{f(0)}{0^2}$ at the origin removes the singularity at the origin, making $\frac{f(x)}{x^2}$ continuous and positive so the extreme value theorem implies it achieves a positive minimum on the interval.

\section{Stochastic Gradient Descent under the PL Condition (Again)}
\label{sec:B}

In the original paper, we give the following result for minimizing a function of the form $f(x) = \mathbb{E}[f_i(x)]$ using stochastic gradient descent (SGD) iterations $x_{k+1}=x_k-\alpha_k\nabla_{i_k}f(x_k)$ for a randomly chosen $i_k$.
\begin{theoremNN}
Assume that $f$ has an $L$-Lipschitz continuous gradient, $f$ has a non-empty solution set $\mathcal{X}^*$, $f$ satisfies the PL inequality (\ref{pl_ineq}), and $\mathbb{E}[\norm{\nabla f_i(x)}^2] \leq C^2$ for all $x$ and  some $C$. If we use the SGD algorithm with $\alpha_k = \frac{2k+1 }{ 2\mu(k+1)^2}$, then we get a convergence rate of 
\ben
\mathbb{E}[f(x_k) - f^*] \leq \frac{L C^2 }{ 2 k\mu^2}.
\een
If instead we use a constant $\alpha_k = \alpha < \frac{1}{2\mu}$, then we obtain a linear convergence rate up to a solution level that is proportional to $\alpha$,
\[
\mathbb{E}[f(x_k) - f^*] \leq (1-2\mu\alpha)^k[f(x_0) - f^*] + \alpha\frac{LC^2}{4\mu}.
\]
\end{theoremNN}
This result shows that under PL SGD achieves the convergence rates of SGD for strongly-convex functions: that with a decreasing step size we obtain an error of $O(1/k)$ and that with a constant step size $\alpha$ we converge linearly to a region of size $O(\alpha)$.
The problem with this result is that the assumptions cannot hold globally. The PL condition requires that the gradient is unbounded but the assumption that $\mathbb{E}[\norm{\nabla f_i(x)}^2] \leq C^2$ requires that the gradient is bounded. Indeed, the latter assumption does not hold for least squares. In the paper we avoided this issue by only assuming that the bound only holds at the particular $x_k$ values. Below we give a variant on the result that instead assumes that the variation in the gradients is bounded in the sense that $\mathbb{E}[\norm{\nabla f_i(x) - \nabla f(x)}^2] \leq \sigma^2$ for some constant $\sigma$ for all $x$.
\begin{theorem}
Assume that $f$ has an $L$-Lipschitz continuous gradient, $f$ has a non-empty solution set $\mathcal{X}^*$, $f$ satisfies the PL inequality (\ref{pl_ineq}), and $\mathbb{E}[\norm{\nabla f_i(x) - \nabla f(x)}^2] \leq \sigma^2$ for all $x$ and  some $\sigma$. If we use the SGD algorithm with $\alpha_k = \frac{2}{\mu(k+\gamma)}$ with $\gamma = \frac{2L}{\mu}$, then we get a convergence rate of 
\[
\mathbb{E}[f(x_k) - f^*] \le \frac{2L\max\left\{ f(x_0) - f^*, \frac{\sigma^2}{\mu}\right\}}{\mu k + 2L}.
\]
If instead we use a constant $\alpha_k = \alpha \leq \frac{1}{L}$, then we obtain a linear convergence rate up to a solution level that is proportional to $\alpha$,
\[
\mathbb{E}[f(x_k) - f^*] \leq (1-\mu\alpha)^k[f(x_0) - f^*] + \alpha\frac{L\sigma^2}{2\mu}.
\]
\end{theorem}
\begin{proof}
We start by following the classic approach of~\citet{ghadimi2013stochastic} to bounding the progress of SGD for $L$-smooth functions, which starts by plugging in the SGD iteration into the descent lemma to give
\[
f(x_{k+1}) \le f(x_k) - \alpha_k \langle \nabla f(x_k), \nabla f_{i_k}(x_k) \rangle + \frac{L\alpha_k^2}{2} \Vert{}\nabla f_{i_k}(x_k)\Vert{}^2.
\]
Taking the conditional expectation with respect to $i_k$ at time $k$ and using $\mathbb{E}[\nabla f_{i_k}(x_k)]=\nabla f(x_k)$ gives
\[
\mathbb{E}_k[f(x_{k+1})] \le f(x_k) - \alpha_k \Vert{}\nabla f(x_k)\Vert{}^2 + \frac{L\alpha_k^2}{2} \mathbb{E}_k\left[\Vert{}\nabla f_{i_k}(x_k)\Vert{}^2\right].
\]
Now substitute the variance-like decomposition
\[
\mathbb{E}_k\left[\Vert{}\nabla f_{i_k}(x_k)\Vert{}^2\right] = \Vert{}\nabla f(x_k)\Vert{}^2 + \mathbb{E}_k\left[\Vert{}\nabla f_{i_k}(x_k) - \nabla f(x_k)\Vert{}^2\right] \le \Vert{}\nabla f(x_k)\Vert{}^2 + \sigma^2,
\]
and group terms to obtain
\[
\mathbb{E}_k[f(x_{k+1})] \le f(x_k) - \alpha_k \left( 1 - \frac{L\alpha_k}{2} \right) \Vert{}\nabla f(x_k)\Vert{}^2 + \frac{L\alpha_k^2\sigma^2}{2}.
\]
Subtract $f^*$ and apply the PL inequality (assuming that $\alpha_k < \frac 2 L$) to obtain the inequality we use to analyze the various step sizes,
\begin{equation}
\mathbb{E}_k[f(x_{k+1}) - f^*] \le \left( 1 - 2\mu\alpha_k \left( 1 - \frac{L\alpha_k}{2} \right) \right) (f(x_k) - f^*) + \frac{L\alpha_k^2\sigma^2}{2},
\label{eq:SGDprog}
\end{equation}
or if $\alpha_k \leq \frac 1 L$ we can use that $\left( 1 - \frac{L\alpha_k}{2} \right) \geq \half$ to get
\begin{equation}
\mathbb{E}_k[f(x_{k+1}) - f^*] \le \left( 1 - \mu\alpha_k\right) (f(x_k) - f^*) + \frac{L\alpha_k^2\sigma^2}{2},
\label{eq:SGDprog2}
\end{equation}
\textbf{Constant step size}. Setting $\alpha_k = \alpha \leq 1/L$ and recursing using that $\sum_{i=0}^{k-1}\rho^i \leq \sum_{i=0}^\infty \rho^i = \frac{1}{1-\rho}$ gives
\begin{align*}
\mathbb{E}[f(x_k) - f^*] & \le \left( 1 - \mu\alpha  \right)^k [f(x_0) - f^*] + \frac{L\alpha^2\sigma^2}{2}\sum_{i=0}^{k-1}\left (1-\mu\alpha\right)^i\\
& \le \left( 1 - \mu\alpha \right)^k [f(x_0) - f^*] + \frac{L\alpha\sigma^2}{2\mu }.
\end{align*}
\textbf{Decreasing step size}. We will set $\alpha_k = \frac{2}{\mu(k+\gamma)}$ with $\gamma = \frac{2L}{\mu}$ which guarantees that $\alpha_k \leq 1/L$. Using this choice in~\eqref{eq:SGDprog2} and simplifying gives
\begin{align*}
\mathbb{E}[f(x_{k+1}) - f^*] & \le \left(1 - \frac{2}{k+\gamma}\right)\mathbb{E}[f(x_k)-f^*] + \gamma\frac{\sigma^2}{\mu(k+\gamma)^2}.
\end{align*}
Define $\zeta = \gamma\max\left\{f(x_0)-f^*, \frac{\sigma^2}{\mu}\right\}$. We will show that $\mathbb{E}[f(x_k)-f^*] \leq \frac{\zeta}{k+\gamma}$ for $k\geq 0$. The base case is that $f(x_0)-f^* \leq \frac{\zeta}{\gamma}$ which holds because it is equivalent to $f(x_0)-f^* \leq \max\left\{f(x_0)-f^*, \frac{\sigma^2}{\mu}\right\}$. For $k\geq 0$ we have
\begin{align*}
\mathbb{E}[f(x_{k+1}) - f^*] & \le \frac{k+\gamma - 2}{k+\gamma}\frac{\zeta}{k+\gamma} + \frac{\zeta}{(k+\gamma)^2}\\
& = \zeta\frac{k + \gamma - 1}{(k + \gamma)^2}.
\end{align*}
We have $\frac{k+\gamma-1}{(k+\gamma)^2} < \frac{1}{k+\gamma+1}$ since $(k+\gamma-1)(k+\gamma+1) = (k+\gamma)^2 - 1< (k+\gamma)^2$. Thus we have that
\begin{align*}
\mathbb{E}[f(x_k) - f^*] & \le \frac{\gamma\max\left\{f(x_0) - f^*, \frac{\sigma^2}{\mu} \right\}}{k + \gamma} \\
& = \frac{2L\max\left\{ f(x_0) - f^*, \frac{\sigma^2}{\mu}\right\}}{\mu k + 2L}.
\end{align*}
\end{proof}
We note that similar results have appeared in the literature~\citep{bottou2018optimization,gower2021sgd,garrigos2023handbook}, but that these assumptions lead to a particularly-clean result that holds under the PL condition. Rather than assuming $\mathbb{E}[\norm{\nabla f_i(x) - \nabla f(x)}^2] \leq \sigma^2$ for all $x$, Robert Gower pointed out to us that it is also possible to obtain a similar result assuming that the noise inequality holds only for the solutions $x^*$ but assuming that the individual functions $f_i$ are themselves $L$-smooth. A similar argument for the case of strongly-convex functions appears in~\citet{needell2014stochastic}.
 We highlight that these results hold for any iteration and do not require averaging of the iterations, though they do require the step sizes to be sufficiently small. On the other hand, to achieve the optimal $O(1/k)$ rate, the decreasing step size sequence above requires knowing both $L$ and $\mu$. We can obtain the $O(1/k)$ rate without a dependence of $\mu$ by using exponential decreasing step sizes~\citep{li2021second}. We can also obtain the $O(1/k)$ rate without a dependence on $L$ and $\mu$ using exponential step sizes~\citep{li2026towards}, although such step sizes may lead to an exponential increase in the function value on early iterations.
\end{document}